\documentclass[12pt,a4paper,leqno,twoside]{article}

\usepackage{amsfonts}
\usepackage{amsmath}
\usepackage{amssymb}
\usepackage{amstext}
\usepackage{amsthm}   % proof environment :)
\newcommand{\R}{\mathbb R}

\newcommand{\W}{{\cal{W}}}
\newcommand{\V}{{\cal{V}}}
\newcommand{\C}{{\cal{C}}}

\newcommand{\A}{{\bf {A}}}

\newcommand{\D}{{\cal{D}}}
\newcommand{\So}{{\cal{S}}}
\newcommand{\F}{{\cal{F}}}
\newcommand{\intrn}{\int_{\R^n}}

\newtheorem{teo}{Theorem}

\newtheorem{lemma}[teo]{Lemma}

\newtheorem{teoA}{Theorem A.$\!\!$}
\newtheorem{lemA}[teoA]{Lemma A.$\!\!$}

\theoremstyle{definition}

\theoremstyle{remark}
\newtheorem{rem}[teo]{Remark}

\title{A multiplicity result for\\
a semilinear Maxwell type equation}
\author{Antonio Azzollini}
\date{Dipartimento di Matematica\\
Universit\`a di Bari,\\
Via Orabona 4, 70125 Bari, Italy\\
e-mail: azzollini@dm.uniba.it}
\usepackage{amsmath}

\begin{document}
\maketitle
%\begin{abstract}
%To put abstract
%\end{abstract}
\begin{abstract}
In this paper we look for solutions of the equation
    \begin{equation*}
        \delta d\A=f'(\langle\A,\A\rangle)\A\hbox{ in }\R^{2k},
    \end{equation*}
where $\A$ is a 1-differential form and $k\geq 2.$ These solutions
are critical points of a functional which is strongly degenerate
because of the presence of the differential operator $\delta d.$
We prove that, assuming a suitable convexity condition on the
nonlinearity, the equation possesses infinitely many finite energy
solutions.
\end{abstract}

\section*{Introduction}\label{sec:intro}
It is well known that the Maxwell equations in the empty space,
written by the differential forms language, are the Euler-Lagrange
equations of
the following action functional%
    \begin{equation}\label{Sforme}%
        \mathcal{S}=\int_{\R^4}\left\langle d\eta,d\eta\right\rangle
        \sigma.
    \end{equation}
Here
    \begin{equation*}
        \eta  = \sum_{i=1}^{3}A_{i}dx^{i}+\varphi dt,\quad
        A_{i},\varphi   : \R^{4}\rightarrow\R,
    \end{equation*}
is the gauge potential 1-form in the space-time $\R^{4}$, $d\eta$
denotes the exterior derivative of $\eta$,  $\sigma$ is the volume
form, and for any differential form $\gamma$
    \begin{equation*}
        \left\langle \gamma,\gamma\right\rangle :=\ast\left(  \ast \gamma\wedge
        \gamma\right)
    \end{equation*}
where $\ast$ is the Hodge operator with respect to the Minkowski
metric in $\R^{4}$.

According to the classical theory of the electrodynamics, when the
electromagnetic field is generated by an assigned source $j$ (e.g.
a particle matter), then the action functional becomes
    \begin{equation*}%
        \mathcal{S}=\int_{\R^4}(\left\langle
        d\eta,d\eta\right\rangle-\langle j,\eta\rangle)
        \sigma.
    \end{equation*}
When instead the source of the field is not assigned but it is an
unknown of the problem, then there are two opposite mathematical
models describing the interaction between the electromagnetic
field and its source: the dualistic model and the unitarian model.

The dualistic model consists in coupling the Maxwell equation with
another field equation describing the dynamics of the source that
is represented by a travelling solitary wave (i.e. a solution of a
field equation whose energy density travels as a localized
packet). This approach has been analyzed in many papers and
several existence and multiplicity results have been obtained (see
e.g. \cite{B.F.4}, \cite{D'AM}, \cite{Coc} and \cite{CocG}).

More recently, an unitarian field theory has been introduced by
Benci and Fortunato \cite{B.F.}, following an idea from Born and
Infeld (see \cite{B.I.}).  According to this theory (we refer to
\cite{B.F.} and \cite{B.F.2} for more details), electromagnetic
field and matter field are both expression of only one physical
entity, and the interaction between them is described by
introducing a nonlinear Poincar\'e invariant perturbation in the
Maxwell Lagrangian in the empty space.

Following this new unitarian theory, we perturb the Lagrangian in
\eqref{Sforme} adding a nonlinear term and obtaining the modified
action functional
    \begin{equation*}\label{Sforme2}
        \mathcal{S}=\int_{\R^4}\Big(\left\langle d\eta,d\eta\right\rangle
        -f(\langle\eta,\eta\rangle)\Big)\,\sigma
    \end{equation*}
where $f:\R\rightarrow\R.$\\
The Euler-Lagrange equation is the following
nonhomogeneous$\phantom{l}$ Maxwell equation
    \begin{equation}\label{one}
        \delta d\eta=j(\eta)
    \end{equation}
where
    \begin{equation}\label{eq:nonhom}
        j(\eta)=f'(\langle\eta,\eta\rangle)\eta
    \end{equation}
and $\delta=\ast d \ast.$\\
The $1-$form $j$ representing the source depends itself on the
gauge $1-$form $\eta$, so the equation \eqref{eq:nonhom} describes
the dynamics of the electromagnetic field in presence of an
auto-induction phenomenon.\\
From now on, we will refer to \eqref{eq:nonhom} as the semilinear
Maxwell equation (SME).

In \cite{ABDF} the equation \eqref{eq:nonhom} has been considered
in the magnetostatic case, namely when it has the form
    \begin{equation}\label{eq:two}
        \delta d\A=f'(\langle\A,\A\rangle)\A
    \end{equation}
where
    \begin{equation*}
        \A   =\sum_{i=1}^{3}A_{i}dx^{i},\quad  A_{i}   : \R^{3}\rightarrow\R,
    \end{equation*}
and the metric on $\R^3$ is the euclidean one. In that paper a
solution $\A$, with the property $\delta\A=0$, has been found. In
\cite{A2}, ignoring the physical origin of the problem, the
equation \eqref{eq:two} has been studied in the more general
context of the $k-$forms on a $n-$Riemannian manifold $M$, and a
multiplicity result has been proved when $M$ is compact.\\
In the
same spirit of that paper, here we consider the problem just from
a mathematical point of view, looking for solution of
    \begin{equation}\label{eq:myprob}
        \left\{
            \begin{array}{l}
                \delta d\A =f'(\langle\A,\A\rangle)\A\\
                \\
                \A   =\sum_{i=1}^{n}A_{i}dx^{i},\quad  A_{i}   :
                \R^{n}\rightarrow\R
            \end{array}
        \right.
    \end{equation}
where we consider $\R^n$ endowed with the euclidean metric. In the
sequel we often will use the notation $\A$ to denote also the
vector field $(A_1,A_2,A_3).$

Now, consider $n\geq1$ even and denote by $\Lambda^1(\R^n)$ the
set of the $1-$forms on $\R^n$ with compact support and by $T$ the
group of transformations on $\R^n$ so defined
    \begin{multline}\label{eq:defT}
                g\in T\Longleftrightarrow g\in O(n) \hbox{ and }\\
                \exists
                (g_i)_{1\leq i\leq n/2}\hbox{ in } O(2) \hbox{ s.t. }
                g=
                \begin{pmatrix}
                    g_1& 0 & \cdots & 0 \\
                    0 & g_2 & \cdots & 0 \\
                    \vdots&\vdots&\ddots&\vdots \\
                    0 & 0 & \cdots & g_{n/2}
                \end{pmatrix}
    \end{multline}
where $O(n)$ and $O(2)$ are respectively the orthogonal groups in
$\R^n$ and $\R^2.$

Moreover, denote by $(\cdot|\cdot)$ the scalar product on $\R^n$
and assume that
    \begin{description}
        \item{$f_1$)} $f\in C^1(\R,\R),$ $f(0)=0$, $\forall t \geq 0: f'(t)\geq
        0$
    \end{description}
and for $2<p<2^*<q$, with $2^*=\frac{2n}{n-2}$,
    \begin{description}
        \item{$f_2$)} $\exists c_1>0$ s.t. $\forall x, y \in
        \R^n$
            \begin{multline}
                f\big((x|x)\big)-f\big((y|y)\big)-
                2 f'\big((y|y)\big)(y|x-y)\\
                \geq
                c_1\min\big((x-y|x-y)^\frac{p}{2},
                (x-y|x-y)^\frac{q}{2}\big),
            \end{multline}
        \item{$f_3$)} $\exists c_2>0$ s.t. $|f'(t)|\leq c_2\min\big(
        t^{\frac{p}{2}-1},t^{\frac{q}{2}-1}\big)$, $\forall t\geq 0$,
        \item{$f_4$)} $\exists R>0$ and $\alpha>2$ s.t.
        $0<\frac{\alpha}{2}\,f(t)\leq f'(t)t,$ $\forall t\geq R.$
    \end{description}

The main result of this paper is the following
    \begin{teo}\label{th:maintheorem}
        Let  $n\geq4$ be even  and assume that $f$
        satisfies ($f_1-f_4$). Then there exist infinitely many
        nontrivial weak solutions of \eqref{eq:myprob}.

        Moreover these
        solutions have the following particular symmetry:
            \begin{equation*}
                \mathbf{A}(x)=g^{-1}\mathbf{A}(gx),\ \ \forall g\in
                T.
            \end{equation*}
    \end{teo}
In the sequel we will assume $f_1,\ldots,f_4$ holding.
    \begin{rem}
        Set $g(x)=f(x^2)$ and suppose $f\in C^2(\R,\R).$
        Observe that by $f_1$ and $f_3$ we deduce $g''(0)=0$,
        and then we find the so called ``zero mass" case. This
        case has been dealt with by Berestycki \& Lions \cite{B.L.1,B.L.2} and
        more recently by Pisani \cite{P}.
    \end{rem}
    \begin{rem}\label{rem:conv}
        For every $x\in \R^n$ we can
        define the scalar product $\langle\cdot,\cdot\rangle_x$ on the
        vector space $\Lambda^1(\R^n).$
        The assumption $f_2$ is a condition on the convexity of the functional
            $$I_x(\xi)=f(\langle\xi,\xi\rangle_x).$$
        In fact, if we take
        $\xi,\psi\in\Lambda^1(\R^n),$ $\lambda\in]0,1[$ and set
        $\eta=\lambda\xi+(1-\lambda)\psi$, by $f_2$ we have
            \begin{multline}
                \lambda \big(f(\langle\xi,\xi\rangle_x) - f(\langle\eta,\eta\rangle_x)-  2
                f'(\langle\eta,\eta\rangle_x)\langle\eta,\xi-\eta\rangle_x\big)\\
                \geq   \lambda\overline
                c\,
                \min\big(\langle\xi-\eta,\xi-\eta\rangle_x^\frac{p}{2},
                \langle\xi-\eta,\xi-\eta\rangle_x^\frac{q}{2}
                \big)>0
                \label{eq:firstconv}
            \end{multline}
        and
            \begin{multline}
                (1-\lambda)\big(f(\langle\psi,\psi\rangle_x)-f(\langle\eta,\eta\rangle_x)
                -
                2f'(\langle\eta,\eta\rangle_x)\langle\eta,\psi-\eta\rangle_x\big)\\
                \geq
                (1-\lambda)\overline
                c\,
                \min\big(\langle\psi-\eta,\psi-\eta\rangle_x^\frac{p}{2},
                \langle\psi-\eta,\psi-\eta\rangle_x^\frac{q}{2}
                \big)>0
                \label{eq:secondconv}.
            \end{multline}
        Since
            $$\lambda f'(\langle\eta,\eta\rangle_x)\langle\eta,\xi-\eta\rangle_x+(1-\lambda)
            f'(\langle\eta,\eta\rangle_x)\langle\eta,\psi-\eta\rangle_x=0,$$
        adding \eqref{eq:firstconv} to \eqref{eq:secondconv} we obtain
            \begin{equation}
                \lambda f(\langle\xi,\xi\rangle_x)+(1-\lambda)f(\langle\psi,\psi\rangle_x)-
                f(\langle\eta,\eta\rangle_x)>0
            \end{equation}
        and then for every $x\in\R^n$ the functional $I_x$ is strictly convex.

        The function $f:\R\rightarrow\R$ s.t.
    \begin{equation*}
        f(x)=\left\{
            \begin{array}{lcl}
                a|x|^\frac{p}{2}+b & \hbox{if} & |x|>1\\
                c|x|^\frac{q}{2}   & \hbox{if} & |x|\leq1
            \end{array}
        \right.
    \end{equation*}
where $2<p<2^*<q$ and $(a,b,c)\in\R^2\times]0,+\infty[$ is any
solution of the system
    \begin{equation*}
        \left\{
            \begin{array}{lcl}
                a+b & = & c\\
                ap  & = & cq
            \end{array}
        \right.
    \end{equation*}
    is an example of function satisfying $f_2$ (see the Appendix for details).
    \end{rem}

    The paper is organized as follows:\\
in section 1, following \cite{B.F.}, we will use a new functional
framework related to the Hodge decomposition of the vector field
$\A.$ We will be led  to study the problem in the space
    \begin{equation*}
        \mathcal{D}(\R^n):=\left\{  u\in L^{6}(\R^n):%
        {\displaystyle\int}
        \left\vert \nabla u\right\vert ^{2}dx<+\infty\right\}
    \end{equation*}
and in the Orlicz space $L^{p}+L^{q}$ ($2<p<6<q$). We will recall
some basic theorems, obtained in \cite{B.F.,P}, describing the
relations between these spaces, and two results, proved
respectively in \cite{P} and \cite{AP}, which will be necessary to
get regularity and compactness.

In section 2, we will give a proof of Theorem
\ref{th:maintheorem}, using a well known multiplicity abstract
result (see \cite{A.R.,B.B.F.}). Assumption $f_2$ will play a key
role in order to get regularity.

Finally, in the appendix we will show an example of function
satisfying the assumptions of Theorem \ref{th:maintheorem}.

\section{The functional setting}

From now on, taken $\A=\sum_{i=1}^{n}A_idx^i$  a $1-$form, by
$\nabla\A$ we mean the Jacobian matrix of the field
$(A_1,A_2,\ldots,A_{n} )$ and if ${\bf B}$ is another $1-$form we
will use the notation $(\nabla\A|\nabla{\bf B})$ to mean the
product
    \[ \left(  \nabla \A\mid\nabla
        {\bf B}\right)  =Tr\left[  \left(  \nabla \A\right) \left(  \nabla
        {\bf B}\right)  ^{T}\right]
    \]
where $\left(  \nabla {\bf B}\right)^{T}$ is transposed of $\nabla
{\bf B}$ and $Tr$ denotes the trace. Moreover in the sequel we
will write $(\A|{\bf B})$ to mean the scalar product between $\A$
and ${\bf B}$ and we will use $|\A|^2$ and $|\nabla\A|^2$ in the
place of $( \A|\A)$ and $(\nabla\A|\nabla\A).$

The functional of the action associated to \eqref{eq:two} is
    \begin{equation}\label{energy}
        J(\A)=\frac{1}{2}\intrn\langle d\A,d\A\rangle\,dx-\frac{1}{2}\intrn
        f(|\A|^2)\,dx
    \end{equation}
where $dx=dx^1\wedge dx^2\wedge\ldots\wedge dx^n,$ being
$\{dx^1,\ldots,dx^n\}$ the canonical basis of $\Lambda^1(\R^n).$
The strongly degenerate nature of the functional $J$ doesn't allow
us to approach this problem in a standard way. In other words, the
functional $J$ doesn't present the geometry of the mountain pass
in any space with finite codimension. This strongly indefiniteness
of the functional depends on the fact that, in general,
    \begin{equation*}
        \intrn\langle d\A,d\A\rangle\,dx\neq \intrn |\nabla
        \A|^2\,dx
    \end{equation*}
since the equality holds only if $\delta \A=0.$ As a consequence,
we
don't have an a priori bound on the norm $\|\nabla \A\|_{L^2}.$\\
To overcome this difficulty, we look to the Hodge decomposition
theorem of the differential forms in order to split
    \begin{equation}\label{split}
        \A=u+d w=u+\nabla w
    \end{equation}
where $u$ is a $1-$form  s.t.
    \begin{equation}\label{defu}
        \delta u=0
    \end{equation}
and $w$ is a $0-$form, i.e. $w:\R^n\to\R.$\\
Substituting the splitting \eqref{split} in
\eqref{energy}, we obtain
    \begin{equation}\label{defJ}
        J(u,w):=J(u+d w)=\frac{1}{2}\intrn |\nabla u|^2\,dx-\frac{1}{2}
        \intrn f(|u+\nabla
        w|^2)\,dx.
    \end{equation}
Now we introduce the spaces where the
functional $J$ is defined.\\
For $2<p<\frac{2n}{n-2}<q$, denote by $(L^p(\R^n),|\cdot|_p)$ and
$(L^q(\R^n),|\cdot|_q)$ the Lebesgue spaces defined as the closure
of $\Lambda^1(\R^n)$ with respect to the norm
    \begin{equation*}
        |\xi|_h=\Big(\intrn
        |\xi|^h\,dx\Big)^\frac{1}{h},\quad h=p,q.
    \end{equation*}
Consider the space
    $$L^{p}+L^{q}:=\{\xi\mid\exists \xi_1\in L^p(\R^n)\hbox{ and }\xi_2\in
    L^q(\R^n) \hbox{ such that }\xi=\xi_{1}+\xi_{2}\}.$$
It is well known that $L^p+L^q$ is a Banach space with the norm
\begin{equation}
\left\|  \xi\right\|  _{L^{p}+L^{q}}=\inf\left\{  \left\|
\xi_{1}\right\| _{L^{p}}+\left\|  \xi_{2}\right\|
_{L^{q}}:(\xi_1,\xi_2)\in L^p\times
L^q,\;\xi_{1}+\xi_{2}=\xi\right\}
\label{normaa}%
\end{equation}
and its dual space is $L^{p'}\cap L^{q'}$, where
$p'=\frac{p}{p-1}$ and $q'=\frac{q}{q-1},$ endowed with the norm
    \begin{equation*}
        \|\xi\|_{L^{p'}\cap L^{q'}}    :=
        \|\xi\|_{L^{p'}}+\|\xi\|_{L^{q'}}.
    \end{equation*}\\
Denote by $C^\infty_0(\R^n)$ the space of the smooth functions
with compact support, and set
    \begin{eqnarray*}
        %L^q_k(M)        &:=&\overline{\Lambda^k(M)}^{|\cdot|_q},\\
        \D(\R^n)        &:=&\overline{\Lambda^1(\R^n)}^{\|\cdot\|},\\
        \D^{p,q}(\R^n)    &:=&\overline{C^\infty_0(\R^n)}^{\|\cdot\|_{p,q}}
    \end{eqnarray*}
where, for every $\xi\in\Lambda^1(\R^n)$,
    \begin{eqnarray*}
        %|\xi|_q^q
        %&:=&\int_M\langle\xi,\xi\rangle^{\frac{q}{2}}\,dx.\nonumber\\
        %\nonumber\\
        \|\xi\|^2       &:=& \int_{\R^n}\langle d \xi,d
        \xi\rangle\,dx
        +\intrn\langle\delta\xi,\delta\xi\rangle\,dx\\
        \noalign{\hbox{and for every  $g\in C^\infty_0(\R^n)$}}
        \|g\|_{\D^{p,q}} &:=& \|\nabla g\|_{L^p+L^q}
    \end{eqnarray*}

We recall some results on the space $L^p+L^q$
    \begin{teo}\label{teo:onLp+Lq}
        \begin{enumerate}
        \item $\Lambda^1(\R^n)$ is dense in $L^p+L^q.$
        \item Let $\xi\in L^p +L^q$ and set
            \begin{equation}\label{defomega}
                \Omega_\xi:=\big\{x
                \in\R^n\big|\,|\xi(x)|>1\big\}.
            \end{equation}
        Then
        \begin{multline}\label{controlnorm}
                \max\left(\|\xi\|_{L^q(\R^n-\Omega_\xi)}-1,
                \frac{1}{1+|\Omega_\xi|^{1/r}}\|\xi\|_{L^p(\Omega_\xi)}\right)\\
                \leq\|\xi\|_{L^p+L^q}\leq\max\left(\|\xi\|_{L^q(\R^n-\Omega_\xi)},
                \|\xi\|_{L^p(\Omega_\xi)}\right)
        \end{multline}
        where $r=p\,q/(q-p).$
        \item For every $r\in [p,q]:$ $L^r\hookrightarrow L^p+L^q$
        continuously.
        \item The embedding
            \begin{equation}\label{emb}
                \D(\R^n)\hookrightarrow L^p+L^q
            \end{equation}
        is continuous.
        \item Set
            \begin{equation}\label{F}
                \F:=\left\{\xi:\R^n\rightarrow\R^n|\forall g\in
                T,
                \hbox{ for a.e. }
                x\in\R^n: \xi(gx)=g\xi(x)\right\}
            \end{equation}
            where $T$ is defined by \eqref{eq:defT}, and define the
            space $\D_r(\R^n)$ as follows
            \begin{equation}\label{Dr}
                \D_r(\R^n):=\D(\R^n)\cap \F.
            \end{equation}
        Then $\D_r(\R^n)\hookrightarrow L^p+L^q$ compactly.
        \end{enumerate}
    \end{teo}
    \begin{proof}
        \begin{enumerate}
        \item It can be easily showed using the definition of the
        $L^p+L^q$-norm and the density of $\Lambda^1(\R^n)$ in the
        spaces $L^p(\R^n)$ and $L^q(\R^n).$
        \item See Lemma 1 in \cite{B.F.}.
        \item See Corollary 9 in \cite{P}.
        \item It follows from 3 and the Sobolev continuous embedding
            $$\D(\R^n)\hookrightarrow L^\frac{2n}{n-2}.$$
        \item The proof follows combining a compactness theorem
        presented in \cite{AP} (see Theorem A.$1$
        in the Appendix) and Lemma 14 in \cite{B.F.}.
        \end{enumerate}
    \end{proof}
For all $\A\in L^p+L^q,$ consider the functional $F$ defined as
follows
    \begin{eqnarray}
        F(\A)&  := &  \int_{\R^n}f(|\A|^2)\,dx.\label{defF}
    \end{eqnarray}

The following results have been proved in \cite{P}
    \begin{teo}\label{diffF}
        If $f_3$ holds, then the functional $F$ is continuously
        differentiable, and its Frechet differential is the
        continuous and bounded map
            \begin{equation}\label{NemOpf'}
                DF:\A\in L^p+L^q\mapsto
                2\intrn f'(|\A|^2)(\A|\cdot)\,dx\;
                \in (L^p+L^q)'.
            \end{equation}
    \end{teo}
        Using the fact that $f(0)=0,$ from $f_2$ we deduce that
        for every $\xi\in L^p+L^q$
            \begin{equation*}
                f(\langle \xi,\xi\rangle)
                \geq
                c_1\min\big(\langle \xi,\xi\rangle^\frac{p}{2},
                \langle \xi,\xi\rangle^\frac{q}{2}\big)
            \end{equation*}
        pointwise almost everywhere in $\R^n.$
        On the other hand, from $f_1$ and $f_3$ it follows that
            \begin{equation*}
                f(\langle \xi,\xi\rangle)
                \leq
                c_2\min\big(\langle \xi,\xi\rangle^\frac{p}{2},
                \langle \xi,\xi\rangle^\frac{q}{2}\big),
            \end{equation*}
        pointwise almost everywhere in $\R^n.$\\
        So, for every $\xi\in L^p+L^q$
            \begin{equation}\label{growthf2}
                c_1\min\big(\langle \xi,\xi\rangle^\frac{p}{2},
                \langle \xi,\xi\rangle^\frac{q}{2}\big)\leq
                f(\langle \xi,\xi\rangle)
                \leq
                c_2\min\big(\langle \xi,\xi\rangle^\frac{p}{2},
                \langle \xi,\xi\rangle^\frac{q}{2}\big),
            \end{equation}
        and then we deduce that for any $\xi\in L^p+L^q:$\
    \begin{multline}\label{growthf}
        c_1\Big(\int_{\Omega_\xi} |\xi|^p\,dx+\int_{\R^n-\Omega_\xi}
        |\xi|^q\,dx\Big)\\
        \leq \intrn f(\xi)\leq c_2
        \Big(\int_{\Omega_\xi} |\xi|^p\,dx+\int_{\R^n-\Omega_\xi}
        |\xi|^q\,dx\Big)
    \end{multline}
By \eqref{controlnorm} and \eqref{emb},
    \begin{equation}\label{Jfinite}
        \forall u\in \D(\R^n),w\in\D^{p,q}(\R^n):
        J(u,w)<+\infty .
    \end{equation}
In order to have compactness for $J,$ we are going to restrict the
domain of the functional to a subspace
$H\subset\D(\R^n)\times\D^{p,q}(\R^n)$ s.t. for all $(u,w)\in H$
we have that $u+\nabla w\in \F.$\\
It is easy to see that, if we set
    \begin{equation}
        \F':=\left\{w:\R^n\rightarrow\R|\forall g\in  T\hbox{ and for
        a.e. }
        x\in\R^n:w(gx)=w(x)\right\},
    \end{equation}
then, for $w:\R^n\rightarrow\R$ sufficiently smooth, we have
    \begin{equation}\label{rad}
         w\in \F'\Longrightarrow \nabla w\in \F.
    \end{equation}
So, taking \eqref{defu} into account, we set
    \begin{equation}\label{defV}
        \V :=\left\{u\in\D_r(\R^n)|\delta
        u=0\right\}
    \end{equation}
and
    \begin{equation}\label{defW}
        \W:=D^{p,q}(\R^n)\cap \F',
    \end{equation}
and we take $H=\V\times\W.$\\
Observe that $H$ is nonempty. In fact, $\W\neq \emptyset$ and for
any $(a_i)_{1\leq i\leq n/2}$ in $C_0^\infty(\R^n)\cap\F'$ the
$1-$form
    \begin{equation*}
        \xi=\sum_{i=1}^{n/2}a_i(x_{2i-1}dx_{2i}-x_{2i}dx_{2i-1})
    \end{equation*}
belongs to $\V.$

Now, for every $u\in\V$ and $w\in\W$,
        set
    \begin{eqnarray}
        F_u  &  :  &  w\in\W\mapsto F(u+\nabla w)\in\R\label{defFu}\\
        F_w  &  :  &  u\in\V\mapsto F(u+\nabla
        w)\in\R\label{defFw}\\
        J_u  &  :  &  w\in\W\mapsto J(u,w)\in\R\label{defJu}\\
        J_w  &  :  &  u\in\V\mapsto J(u,w)\in\R.\label{defJw}
    \end{eqnarray}

    \begin{rem}\label{difffunct}
        Observe that, by Theorem \ref{diffF}, for every $u\in\V$ and $w\in\W$ the
        functionals $J$, $J_u,$  $J_w,$ $F_u$ and $F_w$ are
        continuously differentiable, and the respective Frechet
        differentials are:
        \begin{eqnarray}
            DJ:\V\times\W&\rightarrow& (\V\times\W)'\\
            \nonumber\\
            {DJ_u \atop DF_u}:\W&\rightarrow& \W'\\
            \nonumber\\
            {DJ_w \atop DF_w}:\V&\rightarrow& \V'.
        \end{eqnarray}
        Moreover, if we set
            \begin{eqnarray}
                \frac{\partial J}{\partial w}(u,w)
                &  :=  &  DJ_u(w)\in\W'\label{defpartdervw}\\
                \frac{\partial J}{\partial u}(u,w)
                &  :=  &  DJ_w(u)\in\V',\label{defpartdervu}
            \end{eqnarray}
        by some computations we can see that, for every $u,\overline
        u\in\V$ and $w,\overline w\in\W,$
            \begin{equation}
                \frac{\partial J}{\partial w}(u,w)[\overline
                w]=DJ(u,w)[0,\overline w]    =    -\intrn  f'(|u+\nabla w|^2)
                ( u+\nabla w|\nabla
                \overline w)\,dx,\label{expl1}
            \end{equation}
            \begin{multline}
                \frac{\partial J}{\partial u}(u,w)[\overline
                u]=DJ(u,w)[\overline u,0]    =    \intrn(\nabla u|\nabla\overline
                u)\,dx\\
                -\intrn f'(|u+\nabla w|^2)(u+\nabla w|\overline u)\,dx.\label{expl2}
            \end{multline}
    \end{rem}
Using \eqref{defpartdervw} and \eqref{defpartdervu}, we can show
the variational nature of the problem \eqref{eq:myprob}
    \begin{teo}\label{variatappr}
        If the couple $(u,w)\in\V\times\W$ solves the system
            \begin{eqnarray}
                \frac{\partial J}{\partial w}(u,w)  &  =  &  0\label{sys1}\\
                \frac{\partial J}{\partial v}(u,w)  &  =  &
                0\label{sys2}
            \end{eqnarray}
        then $\A=u+\nabla w\in\F$ is a finite energy, weak solution of
        \eqref{eq:myprob}.
    \end{teo}
    \begin{proof}
        Let $(u,w)\in\V\times\W$ be a solution of \eqref{sys1} and
        \eqref{sys2}. Then, by \eqref{expl1} and \eqref{expl2},
        for any $\overline u\in\V$ and $\overline w\in\W$
            \begin{equation}\label{expl3}
                \intrn  f'(|u+\nabla w|^2)\big(u+\nabla w|\nabla
                \overline w\big)\,dx  =  0,
            \end{equation}
            \begin{equation}\label{expl4}
                \intrn(\nabla u|\nabla\overline
                u)\,dx-\intrn f'(|u+\nabla w|^2)\big(u+\nabla w|\overline u\big)\,dx
                 =  0.
            \end{equation}
        We want to show that
            \begin{equation}\label{decA}
                \A=u+\nabla w
            \end{equation}
        is a weak solution of
        \eqref{eq:myprob}, namely for all $\varphi\in
                 \Lambda^1(\R^n)$
            \begin{equation}\label{wsol}
                DJ(\A)[\varphi]=\intrn\langle d\A,d\varphi\rangle\,dx-\intrn
                 f'(|\A|^2)(\A|\varphi)\,dx=0.
            \end{equation}
        Actually, it is enough to prove \eqref{wsol} just for
        every $\varphi \in \D_r,$ since $\D_r$ is a natural
        constraint for $J.$  In fact observe that, if we denote by $\mathcal T$
        the group
        of isometric transformations on $\D(\R^n)$ defined as follows
            \begin{equation*}
                G\in{\mathcal T}\Longleftrightarrow
                \exists g\in T\hbox{ s.t. }{\substack{
                \forall\A\in\D \\ \hbox{\small{\small {\small for a.e. }}}x\in\R^n}}:G(\A)(x)=g^{-1}\A(gx),
            \end{equation*}
        then $\D_r$ is the subspace of the fix points of $\D(\R^n)$
        under the action of $\mathcal T$ and
            \begin{equation*}
                \forall G\in{\mathcal
                T},\:\forall\A\in\D(\R^n):J(G(\A))=J(\A).
            \end{equation*}
        Then, by the Palais' Principle of symmetric
        criticality (see \cite{Pa}), $\D_r$ is a natural constraint.
        Let $\varphi\in \D_r.$ As in \eqref{split}, we can split the
        function $\varphi$ and obtain
            \begin{equation}\label{split2}
                \varphi=v+d h=v+\nabla h
            \end{equation}
        where $v\in\V$ and $h\in\W$. Writing \eqref{expl3} and
        \eqref{expl4} with respectively  $\overline u=v$ and $\overline w=h$, we get
        \begin{equation}\label{expl5}
            \intrn  f'(|u+\nabla w|^2)\big(u+\nabla w|\nabla h\big)\,dx  =  0,
        \end{equation}
        \begin{equation}\label{expl6}
            \intrn(\nabla u|\nabla v)\,dx-\intrn  f'(|u+\nabla w|^2)\big(u+\nabla w|v\big)\,dx
            =  0,
        \end{equation}
        so, subtracting \eqref{expl5} from \eqref{expl6}, by \eqref{split2}
        we have
        \begin{equation}\label{explsys}
            \intrn\big(\nabla u|\nabla v\big)\,dx-\intrn
             f'(|u+\nabla
            w|^2)\big(u+\nabla
            w|\varphi\big)\,dx=0.
        \end{equation}
        Since $\delta v =0$, then
        \begin{equation}\label{explcurlcurl}
            \delta d\varphi=
            \delta d (v+d h)=
            \delta d v=-\Delta v,
        \end{equation}
        where $-\Delta:=d\delta+\delta d$ is the Laplace-Beltrami
        operator.
        From \eqref{explsys} and \eqref{explcurlcurl}, we
        deduce that
        \begin{align}\label{oneoftheeq}
            \intrn\langle d u,&d\varphi\rangle\,dx
              -
            \intrn  f'(|u+\nabla w|^2)\big(u+\nabla w|\varphi\big)\,dx\nonumber\\
            &  =
            \intrn\langle u,\delta d\varphi\rangle\,dx
            -\intrn f'(|u+\nabla w|^2)\big(u+\nabla w|\varphi\big)\,dx\nonumber\\
            &  =
            -\intrn\big(u|\Delta v\big)\,dx -
            \intrn f'(|u+\nabla w|^2)\big(u+\nabla w|\varphi\big)\,dx\nonumber\\
            &  =
            \intrn\big(\nabla u|\nabla v\big)\,dx-\intrn
            f'(|u+\nabla w|^2)\big(u+\nabla w|\varphi\big)\,dx=0.
        \end{align}

        Since $u\in\D(\R^n)$ and $w\in\D^{p,q}(\R^n),$
        by \eqref{Jfinite} the energy of $\A$ is finite.

        Finally, since $u,\nabla w\in\F,$
        also $\A\in\F.$
    \end{proof}
\section{Proof of the main theorem} Set
    \begin{eqnarray}
        \C_1  &  :=  &  \left\{(u,w)\in\V\times\W\left|
        \frac{\partial J}{\partial w}(u,w)  =  0\right.\right\}\\
        \C_2  &  :=  &  \left\{(u,w)\in\V\times\W\left|
        \frac{\partial J}{\partial u}(u,w)  =  0\right.\right\}.\label{C2}
    \end{eqnarray}
By Theorem \ref{variatappr}, we are interested in finding the
couples $(u,w)\in\C_1\cap\C_2.$\\
Rendering \eqref{C2} explicit  we have that
    \begin{multline}\label{charC1eq}
        (u,w)\in\C_2\Longleftrightarrow\forall \overline u\in\V:
         \intrn(\nabla u|\nabla\overline
         u)\,dx\\
         -\intrn  f'(|u+\nabla w|^2)\big(u+\nabla w|\overline u\big)\,dx
         =  0.
    \end{multline}
The following theorem characterizes the set $\C_1$
    \begin{teo}\label{charC1}
        There exists a compact map
        $\Phi:\V\rightarrow\W$ s.t.
            \begin{equation}\label{intrPhi}
                \C_1=\left\{\big(u,\Phi(u)\big)\left|u\in\V\right.\right\}.
            \end{equation}
        Moreover the map $\Phi$ is characterized by the following
        property:\\
            \begin{equation}
                \begin{array}{c}\label{defPhi}
                    \hbox{for every $u\in\V,$ $\Phi(u)$ is the unique function in
                    $\W$ s.t.}\\
                    \\
                    F_u\big(\Phi(u)\big)=\displaystyle{\min_{w\in\W}F_u(w).}
                \end{array}
            \end{equation}
    \end{teo}
Before we prove the Theorem \ref{charC1}, we need the following
    \begin{lemma}\label{strconvcond}
        If
            \begin{equation}
                \zeta_n   \rightharpoonup   \zeta\hbox{ in
                }L^p+L^q\label{firstcon}
            \end{equation}
        and
            \begin{equation}
                F(\zeta_n)  \rightarrow
                F(\zeta),\label{seccon}
            \end{equation}
        then
            \begin{equation}\label{strcon}
                \zeta_n\rightarrow\zeta\hbox{ in
                }L^p+L^q.
            \end{equation}
    \end{lemma}
    \begin{proof}
        Let $(\zeta_n)_n$ be a sequence in $L^p+L^q$ and $\zeta\in L^p+L^q$ s.t.
        \eqref{firstcon} and \eqref{seccon} hold. Using $\tilde
        f_2$ for $(\zeta_n)_x$ and $(\zeta)_x$ for all $x\in\R^n$ and $n\geq 1,$
        we have that the following inequality
        holds pointwise:
            \begin{multline}\label{f2applAnA}
                f\big(|\zeta_n|^2\big)    -    f\big(|\zeta|^2\big)-2
                f'\big(|\zeta|^2\big)\Big(\zeta\Big|\zeta_n-\zeta\Big)\\
                 \geq   c_1\,\min\left(\left|\zeta_n-\zeta\right|^p,
                \left|\zeta_n-\zeta\right|^q\right).
            \end{multline}
        Set
            \begin{equation*}
                \Omega_n:\left\{x\in\R^n\big|\left|\zeta_n-\zeta\right|>1
                \right\}.
            \end{equation*}
        Integrating in inequality \eqref{f2applAnA}, by Theorem \ref{diffF} we get
            \begin{align}\label{thcon}
                    F(|\zeta_n|^2)  &  -
                    F(|\zeta|^2)- DF\big(\zeta\big)(\zeta_n-\zeta)\nonumber\\
                    &  \geq
                    c_1\,\displaystyle \int_{\Omega_n}\big|\zeta_n-\zeta\big|^p\,dx+
                    c_1\,\int_{\R^n-\Omega_n}\big|\zeta_n-\zeta\big|^q\,dx\nonumber\\
                    &  =   \displaystyle
                    c_1\,\left(\|\zeta_n-\zeta\|_{L^p(\Omega_n)}^p+
                    \|\zeta_n-\zeta\|_{L^q(\R^n-\Omega_n)}^q\right),
            \end{align}
        By \eqref{firstcon}, \eqref{seccon} and \eqref{thcon} we
        have that
            \begin{equation*}
                \|\zeta_n-\zeta\|_{L^p(\Omega_n)}^p+
                    \|\zeta_n-\zeta\|_{L^q(\R^n-\Omega_n)}^q\longrightarrow0,
            \end{equation*}
        and then we get \eqref{strcon} by \eqref{controlnorm}.
    \end{proof}

    \begin{proof}[Proof of Theorem \ref{charC1}]
        Let $u\in\V$ and consider $F_u$ defined as in
        \eqref{defFu}. By Remark \ref{difffunct} and Remark
        \ref{rem:conv},
        $F_u$ is continuous and strictly convex. Then $F_u$ is
        weakly lower semicontinuous.\\
        Moreover $F_u$ is also coercive. In fact, if $w\in\W$ and
        we set
            \begin{equation*}
                \Omega:=\left\{x\in\R^n\big||u(x)+\nabla
                w(x)|>1\right\},
            \end{equation*}
        then, by \eqref{growthf}, we have
            \begin{eqnarray}\label{coerF}
                F_u(w)  &  =  &  \intrn f(|u+\nabla w|^2)\,dx\nonumber\\
                        &  =  &  \int_{\R^n-\Omega}f(|u+\nabla
                        w|^2)\,dx+\int_\Omega f(|u+\nabla
                        w|^2)\,dx\nonumber\\
                        &\geq &  c_1 \int_{\R^n-\Omega}|u+\nabla
                        w|^q\,dx+c_1 \int_\Omega |u+\nabla
                        w|^p\,dx.
            \end{eqnarray}
        By \eqref{coerF} and \eqref{controlnorm} we deduce that
        $F_u$ is coercive and then, by Weierstrass theorem, $F_u$ possesses
        a minimizer in $\W.$\\
        So, let $\Phi$ be the map defined as follows
            \begin{equation}\label{eq:minimizermap}
                \Phi: \V\rightarrow\W \hbox{ s.t. }\forall u\in\V:
                \Phi(u)\hbox{ minimizes }
                F_u.
            \end{equation}
        Since $F_u$ is strictly convex, for all $u\in\V$ the minimizer
        of the functional $F_u$ is unique, and then the map $\Phi$ is well
        defined and satisfies \eqref{defPhi}.

        Now, before we prove the compactness  of $\Phi:\V\rightarrow\W$, first we
        show that the functional
            \begin{equation}\label{newfunc}
                  u\in\V\mapsto\int_{\R^n}f\big(|u+\nabla\Phi(u)|^2\big)\,dx
            \end{equation}
        is weakly continuous.

        Let
            \begin{equation}\label{hypwcont1}
                u_n\rightharpoonup u\hbox{ in }\V,
            \end{equation}
        then, by 5 of Theorem \ref{teo:onLp+Lq},
            \begin{equation}\label{hypscont}
                u_n\rightarrow u\hbox{ in }L^p+L^q.
            \end{equation}
        Since
            \begin{equation*}
                0\leq F\big(u_n+\nabla\Phi(u_n)\big)=
                F_{u_n}\big(\Phi(u_n)\big)\leq F_{u_n}(0)=
                F(u_n),
            \end{equation*}
        by \eqref{hypscont} and the continuity of $F,$ the
        sequence $\left\{F\big(u_n+\nabla\Phi(u_n)\big)\right\}$
        is bounded.
        Since $F$ is coercive, then
            \begin{equation}\label{unbound}
                u_n+\nabla\Phi(u_n)\hbox{ is
                bounded in }L^p+L^q,
            \end{equation}
        so, by \eqref{hypscont},
            \begin{equation}\label{Phiunbound}
                \nabla\Phi(u_n)\hbox{ is
                bounded in }L^p+L^q.
            \end{equation}
        Set
            \begin{eqnarray*}
                \alpha_n  &  :=  &  \intrn
                f\big(|u_n+\nabla\Phi(u_n)|^2\big)\,dx-\intrn
                f\big(|u+\nabla\Phi(u_n)|^2\big)\,dx\label{defalpha}\\
                \beta_n   &  :=  &  \intrn
                f\big(|u_n+\nabla\Phi(u_n)|^2\big)\,dx-\intrn
                f\big(|u+\nabla\Phi(u)|^2\big)\,dx\label{defbeta}\\
                \gamma_n  &  :=  &  \intrn
                f\big(|u_n+\nabla\Phi(u)|^2\big)\,dx-\intrn
                f\big(|u+\nabla\Phi(u)|^2\big)\,dx\label{defgamma}.
            \end{eqnarray*}
        By \eqref{eq:minimizermap}, certainly we have
            \begin{equation}\label{convsucc}
                \alpha_n\leq\beta_n\leq\gamma_n.
            \end{equation}
        Moreover, by Lagrange theorem,
            \begin{eqnarray}
                \alpha_n  &  =  &  \intrn
                \Big(f\big(|u_n+\nabla\Phi(u_n)|^2\big)-
                f\big(|u+\nabla\Phi(u_n)|^2\big)\Big)\,dx\nonumber\\
                &  =  &
                2\intrn f'\big(|\theta_n|^2\big)
                \left(\theta_n\big|u_n-u\right)\,dx\label{Lagr}
            \end{eqnarray}
        where $\theta_n$ is a suitable convex combination of
        $u_n+\nabla\Phi(u_n)$ and $u+\nabla\Phi(u_n).$
        Since $\{u_n\}$ and $\left\{\nabla\Phi(u_n)\right\}$ are bounded
        in $L^p+L^q$, certainly also $\{\theta_n\}$ is bounded in
        $L^p+L^q$. Then, by Theorem \ref{diffF} and
        \eqref{hypscont}, from \eqref{Lagr} we deduce that
            \begin{equation}\label{alphato0}
                \alpha_n\longrightarrow0.
            \end{equation}
        Analogously we also have that
            \begin{equation}\label{gammato0}
                \gamma_n\longrightarrow0,
            \end{equation}
        so, by \eqref{convsucc}, \eqref{alphato0} and \eqref{gammato0}, we get
            \begin{equation*}
                \beta_n\longrightarrow0
            \end{equation*}
        and then \eqref{newfunc} is weakly continuous.\\
        Now, we prove the compactness of $\Phi.$
        Consider again $(u_n)_{n\geq1}$ in $\V$ s.t.
        \eqref{hypwcont1} holds. By
        \eqref{Phiunbound}, there exists
        $w\in\W$ s.t. (up to a subsequence)
            \begin{equation}\label{hypwcont2}
                \nabla \Phi(u_n)\rightharpoonup\nabla w\hbox{ in
                }L^p+L^q.
            \end{equation}
        From \eqref{hypscont} and \eqref{hypwcont2} we deduce that
            \begin{equation}\label{hypwcont3}
                u_n+\nabla\Phi(u_n)\rightharpoonup u+\nabla w\hbox{ in
                }L^p+L^q
            \end{equation}
        so, using the weak continuity of \eqref{newfunc} and the
        weak lower semicontinuity of $F$ we have
            \begin{equation}\label{wminimize}
                \begin{array}{lcl}
                    F_u\big(\Phi(u)\big)  &  =  &F\big(u+\nabla\Phi(u)\big)   =
                    \displaystyle\lim_n F\big(u_n+\nabla\Phi(u_n)\big)\\
                    &  \geq  & F(u+\nabla
                    w)=F_u(w).
                \end{array}
            \end{equation}
        By the uniqueness of the minimizer of $F_u$, from \eqref{wminimize}
        we deduce that $w=\Phi(u)$, so, by
        \eqref{hypwcont3}, we have
            \begin{equation}\label{hypwcont4}
                u_n+\nabla \Phi(u_n)\rightharpoonup u+\nabla \Phi(u)\hbox{ in
                }L^p+L^q.
            \end{equation}
        But using the weak continuity of \eqref{newfunc}, by \eqref{hypwcont1}
        we also have
            \begin{equation}\label{convoffunc}
                \intrn
                f\big(|u_n+\nabla\Phi(u_n)|^2\big)\,dx\longrightarrow\intrn
                f\big(|u+\nabla\Phi(u)|^2\big)\,dx
            \end{equation}
        so, by Lemma \ref{strconvcond}, from \eqref{hypwcont4} and \eqref{convoffunc}
        we deduce that
            \begin{equation}\label{strongconv}
                u_n+\nabla\Phi(u_n)\longrightarrow u+\nabla\Phi(u)\hbox{ in
                }L^p+L^q.
            \end{equation}
        Now, comparing \eqref{strongconv} with \eqref{hypscont},
        we deduce that
            \begin{equation*}
                \Phi(u_n)\longrightarrow\Phi(u)\hbox{ in
                }\W
            \end{equation*}
        and then $\Phi$ is compact.

        Finally, we prove \eqref{intrPhi}. Observe that, since
        $\frac{\partial J}{\partial w}(u,w)=DF_u (w)$, then
            \begin{equation}\label{critpoindFu}
                (u,w)\in\C_1\Longleftrightarrow DF_u (w)  = 0.
            \end{equation}
        But since $F_u$ is convex, its critical points are
        minimizers, and then
            \begin{equation}\label{charccritpoint}
                DF_u (w)  = 0\Longleftrightarrow w=\Phi(u),
            \end{equation}
        so we have \eqref{intrPhi} by \eqref{critpoindFu} and
        \eqref{charccritpoint}.
    \end{proof}
    \medskip
    Consider the functional $\widehat J:\V\rightarrow\R$
        \begin{equation}\label{defhatJ}
            \widehat J(u):=J\big(u,\Phi(u)\big)=\frac{1}{2}\intrn |\nabla u|^2\,dx
            -\frac{1}{2}F\big(u+\nabla\Phi(u)\big).
        \end{equation}
    The following regularity result holds
        \begin{teo}\label{regJ}
            The functional $\widehat J$ is continuously
            differentiable and its Frechet differential
            $D\widehat J:\V\rightarrow\V'$ has this expression
                \begin{equation}\label{exprdJ}
                    D\widehat J(u)[\overline u]=
                    \intrn(\nabla u|\nabla\overline
                    u)\,dx-\intrn  f'(|u+\nabla \Phi(u)|^2)\big(u+\nabla \Phi(u)|\overline
                    u\big)\,dx.
                \end{equation}
        \end{teo}
        \begin{proof}
            Set
                \begin{equation}
                    \widehat F:u\in\V\mapsto F\big(u+\nabla\Phi(u)\big).
                \end{equation}
            We will prove that $\widehat F\in
            C^1$ so that, clearly, also  $\widehat J\in C^1.$

            Let $u\in\V$. We claim that for all $\overline
            u\in\V-\left\{0\right\}$ the functional $\widehat F$ is
            derivable at $u$ in the direction $\overline u,$ and
            the directional derivative (i.e. the G$\hat{\hbox{a}}$teaux derivative
            $D_G\widehat F$) is
                \begin{equation}\label{exdirderiv}
                    D_G\widehat F(u)[\overline u]=2\intrn  f'(|u+\nabla \Phi(u)|^2)
                    \big(u+\nabla \Phi(u)|\overline
                    u\big)\,dx.
                \end{equation}
            In fact, let $t\in\R-\{0\}$ and set
                \begin{eqnarray*}
                    \alpha(t)  &  :=  &  F\big(u+t\overline
                    u+\nabla\Phi(u+t\overline
                    u)\big)-F\big(u+\nabla\Phi(u+t\overline u)\big),\\
                    \beta(t)   &  :=  &  F\big(u+t\overline
                    u+\nabla\Phi(u+t\overline
                    u)\big)-F\big(u+\nabla\Phi(u)\big),\\
                    \gamma(t)  &  :=  &  F\big(u+t\overline
                    u+\nabla\Phi(u)\big)-F\big(u+\nabla\Phi(u)\big).
                \end{eqnarray*}
            By \eqref{defPhi} we know that
                \begin{eqnarray*}
                    F\big(u+t\overline u+\nabla\Phi(u+t\overline u)\big)
                    &  \leq  &  F\big(u+t\overline u+\nabla\Phi(u)\big)\\
                    \noalign{\hbox{and}}
                    F\big(u+\nabla\Phi(u)\big)
                    &  \leq  &  F\big(u+\nabla\Phi(u+t\overline
                    u)\big),
                \end{eqnarray*}
            and then, certainly, for every $t\in\R-\{0\}$
                \begin{equation}\label{alleqbetleqgam}
                    \alpha(t)\leq\beta(t)\leq\gamma(t).
                \end{equation}
            Now, for every $t\in\R-\{0\}$, set
                \begin{eqnarray*}
                    \tilde\alpha(t)  &  =  &
                    \frac{\alpha(t)}{t},\\
                    \tilde\beta(t)  &  =  &
                    \frac{\beta(t)}{t},\\
                    \tilde\gamma(t)  &  =  &
                    \frac{\gamma(t)}{t}
                \end{eqnarray*}
            and observe that \eqref{exdirderiv} means that
                \begin{equation}\label{convtilbeta}
                    \lim_{t\rightarrow0}\tilde\beta(t)=
                    2
                    \intrn  f'(|u+\nabla \Phi(u)|^2)(u+\nabla \Phi(u)
                    |\overline
                    u)\,dx.
                \end{equation}
            From \eqref{alleqbetleqgam} we deduce that
                \begin{eqnarray*}
                    \tilde\alpha(t)\leq\tilde\beta(t)\leq\tilde\gamma(t)
                    &  \hbox{if}  &  t>0,\\
                    \tilde\gamma(t)\leq\tilde\beta(t)\leq\tilde\alpha(t)
                    &  \hbox{if}  &  t<0,
                \end{eqnarray*}
            and then
                \begin{equation}\label{betwal&gam}
                    \min\big(\tilde\alpha(t),\tilde\gamma(t)\big)
                    \leq\tilde\beta(t)
                    \leq\max\big(\tilde\alpha(t),\tilde\gamma(t)\big).
                \end{equation}
            Now, by Lagrange theorem, we have that
                \begin{eqnarray}
                    \tilde\alpha(t)  &  =  &  \frac{\displaystyle
                    \intrn\Big(f\big(|u+t\overline
                    u+\nabla\Phi(u+t\overline
                    u)|^2\big)-f\big(|u+\nabla\Phi(u+t\overline
                    u)|^2\big)\Big)\,dx}{t}\nonumber\\
                    &  =  &  \frac{2\displaystyle\intrn
                    f'\big(|\theta_t|^2\big)\big(\theta_t\big|t\overline
                    u\big)\,dx}{t}=2\intrn f'\big(|\theta_t|^2\big)
                    \big(\theta_t\big|\overline
                    u\big)\,dx\nonumber\\
                    &  =  & DF(\theta_t)[\overline
                    u]\label{convtilalp}
                \end{eqnarray}
            where $\theta_t$ is a suitable convex combination of
            $u+t\overline u+\nabla\Phi(u+t\overline u)$
            and $u+\nabla\Phi(u+t\overline u).$\\
            Since $\Phi$ is  continuous, we have that
                \begin{eqnarray*}
                    \lim_{t\rightarrow0}u+t\overline
                    u+\nabla\Phi(u+t\overline
                    u)&=&u+\nabla\Phi(u)\hbox{ in }L^p+L^q\\
                    \noalign{\hbox{and}}
                    \lim_{t\rightarrow0}u+\nabla\Phi(u+t\overline
                    u)&=&u+\nabla\Phi(u)\hbox{ in }L^p+L^q,
                \end{eqnarray*}
            and then
                \begin{equation}\label{convtheta}
                    \lim_{t\rightarrow0}\theta_t=u+\nabla\Phi(u)
                    \hbox{ in }L^p+L^q.
                \end{equation}
            By continuity, from \eqref{convtilalp} and
            \eqref{convtheta} we deduce that
                \begin{multline}\label{convtilalpha}
                    \lim_{t\rightarrow0}\tilde\alpha(t)=
                    DF(u+\nabla\Phi(u))[\overline u]\\=2
                    \intrn  f'(|u+\nabla \Phi(u)|^2)(u+\nabla \Phi(u)
                    |\overline
                    u)\,dx.
                \end{multline}
            By the same arguments, we can see that
                \begin{equation}\label{convtilgam}
                    \lim_{t\rightarrow0}\tilde\gamma(t)=
                    2
                    \intrn  f'(|u+\nabla \Phi(u)|^2)(u+\nabla \Phi(u)
                    |\overline
                    u)\,dx,
                \end{equation}
            so, by \eqref{betwal&gam}, \eqref{convtilalpha} and \eqref{convtilgam}
             we get \eqref{convtilbeta}, i.e. and the existence of
            the directional derivative.

            Now observe that from \eqref{exdirderiv} we have
                \begin{equation}
                    D_G\widehat F(u)\in\V',\quad\forall u\in\V
                \end{equation}
            and the map
                \begin{equation}\label{mapdF}
                    D\widehat F_G:u\in\V\mapsto2
                    \intrn  f'(|u+\nabla \Phi(u)|^2)(u+\nabla \Phi(u)
                    |\cdot)\,dx\in\V'
                \end{equation}
            is continuous by Theorem \ref{diffF} and the continuity of
            $\Phi$. Then $\widehat F$ is Frechet differentiable,
            and, for all $u,\overline
            u\in\V$
                \begin{equation}\label{exprdF}
                    D\widehat F(u)[\overline u]=
                    2
                    \intrn  f'(|u+\nabla \Phi(u)|^2)(u+\nabla \Phi(u)
                    |\overline
                    u)\,dx.
                \end{equation}
            From \eqref{exprdF} we have \eqref{exprdJ}.
        \end{proof}
        \begin{teo}\label{critJ}
            If $u\in\V$ is a nontrivial critical point of $\widehat J$, then
            $\A=u+\nabla\Phi(u)\in\F$ is a finite energy,
            nontrivial weak solution of \eqref{eq:myprob}.
        \end{teo}
        \begin{proof}
            Let $u\in\V$ be a critical point of $\widehat J$.
            By \eqref{exprdJ} we have that
                \begin{equation}
                    \intrn(\nabla u|\nabla\overline
                    u)\,dx-\intrn  f'(|u+\nabla \Phi(u)|^2)(u+\nabla \Phi(u)|\overline
                    u)\,dx=0
                \end{equation}
            so, by \eqref{charC1eq}, the couple
            $(u,\Phi(u))\in\C_2.$ Since by Theorem
            \ref{charC1} we also have that
            $(u,\Phi(u))\in\C_1,$ then, by Theorem
            \ref{variatappr}, $\A=u+\nabla\Phi(u)$ is a finite
            energy, weak solution.\\
            Moreover, if $u\neq0$, then
                \begin{equation}\label{nontrsolution}
                    u+\nabla\Phi(u)\neq0.
                \end{equation}
            In fact, if
                \begin{equation}\label{contrdhyp}
                    u=-\nabla\Phi(u),
                \end{equation}
            then
                \begin{equation*}
                    -\Delta\Phi(u)=\nabla\cdot u=0,
                \end{equation*}
            and this should imply
                \begin{equation*}
                    \intrn|\nabla\Phi(u)|^2\,dx=0,
                \end{equation*}
            that is
                \begin{equation}\label{contr}
                    \nabla\Phi(u)=0.
                \end{equation}
            But \eqref{contrdhyp} and \eqref{contr} contradict the
            fact that
            $u\neq0$, so \eqref{nontrsolution} holds.
        \end{proof}
By Theorem \ref{critJ} we are reduced to find the critical points
of $\widehat J$, so we are going to study the geometry and the
compactness properties of the functional in order to apply the
symmetrical mountain pass theorem (see \cite{A.R.,B.B.F.}).
    \begin{teo}\label{JsatP.S:}
        $\widehat J$ satisfies the following Palais-Smale (P-S)
        condition:

        if $\{u_n\}\in\V$ is a sequence s.t. for $M\geq0$
            \begin{align}
                \widehat J(u_n)  &  \leq   M,\quad\forall
                n\geq1\label{firP.S.}\\
                \noalign{\hbox{and}}
                D\widehat J(u_n)&\longrightarrow0\label{secP.S.},
            \end{align}
        then $\{u_n\}\in\V$ is precompact.
    \end{teo}
    \begin{proof}
        Let $\{u_n\}\in\V$ be a sequence s.t. \eqref{firP.S.} and
        \eqref{secP.S.} hold. Since $\Phi$ is compact and the embedding
        $\V\hookrightarrow L^p+L^q$ is compact, we have that the
        map \eqref{mapdF} is compact, so, by standard
        arguments we are reduced to prove  that $\{u_n\}$ is
        bounded.

        Rendering \eqref{firP.S.} explicit, we have
            \begin{equation}\label{Jbound}
                \frac{1}{2}\intrn|\nabla u_n|^2\,dx-\frac{1}{2}\intrn
                f(|u_n+\nabla\Phi(u_n)|^2)\,dx\leq M.
            \end{equation}
        Moreover, from \eqref{secP.S.} we deduce that
            $$ D\widehat
            J(u_n)\big[u_n/\|u_n\|_{\D}\big]\longrightarrow0
            $$
        that is there exists $\varepsilon_n\rightarrow0$ s.t.
            \begin{multline}\label{dJgoes0}
                \intrn|\nabla u_n|^2\,dx-\\
                \intrn
                f'(|u_n+\nabla\Phi(u_n)|^2)\big(u_n+\nabla\Phi(u_n)\big|u_n\big)\,dx=
                \varepsilon_n\|u_n\|_{\D}.
            \end{multline}
        Now, by \eqref{defPhi}, certainly we have that for every
        $w\in\W$
            \begin{equation*}
                0=DF_{u_n}(\Phi(u_n))[w]
                =
                \intrn f'(|u_n+\nabla\Phi(u_n)|^2)\big(u_n+\nabla\Phi(u_n)\big|\nabla
                w\big)\,dx,
            \end{equation*}
        so \eqref{dJgoes0} can be written as follows
            \begin{equation}\label{dJgoes02}
                \intrn|\nabla u_n|^2\,dx-\intrn
                f'(|v_n|^2)|v_n|^2\,dx=
                \varepsilon_n\|u_n\|_{\D},
            \end{equation}
        where we have set $v_n=u_n+\nabla\Phi(u_n).$
        Now, multiplying \eqref{Jbound} by $\alpha$ and
        subtracting \eqref{dJgoes02} we get
            \begin{eqnarray}\label{combined}
                \Big(\frac{\alpha}{2}-1\Big)\intrn
                |\nabla u_n|^2\,dx  +
                \intrn\big[f'(|v_n|^2)|v_n|^2-
                \frac{\alpha}{2} f(|v_n|^2)\big]\,dx
                \phantom{\varepsilon_n\|u_n\|_{\D}}\nonumber\\
                \phantom{\Big(\frac{\alpha}{2}-1\Big)\intrn
                |\nabla u_n|^2\,dx  +
                \intrn\big[f'(|v_n|^2)|v_n|^2-
                \frac{\alpha}{2} }
                \leq M-\varepsilon_n\|u_n\|_{\D}.
            \end{eqnarray}
        Using $f_4$,
        from \eqref{combined} we deduce that $\{u_n\}$ is bounded.
        \end{proof}
        \begin{teo}\label{firstgeomount}
            There exist $\rho>0$ and $C>0$ s.t.
                \begin{equation*}
                    \widehat J(u)>C,\quad\forall u\in\V\cap S_\rho,
                \end{equation*}
            where $S_\rho:=\big\{u\in\D\big|\|u\|_{\D}=\rho\big\}.$
        \end{teo}
        \begin{proof}
            Let $u\in\V$ and consider $\Omega_u$ defined as in
            \eqref{defomega}.
            Since $p<2^*<q$ we have that
                \begin{eqnarray}
                    |u(x)|^p  &  \leq  &  |u(x)|^{2^*},\quad\hbox{if
                    }x\in\Omega_u\label{p<2^*}\\
                    |u(x)|^q  &  \leq  &  |u(x)|^{2^*},\quad\hbox{if
                    }x\in\R^n-\Omega_u,\label{q>2^*}
                \end{eqnarray}
            so, by \eqref{p<2^*} and \eqref{q>2^*}, using
            \eqref{growthf}, \eqref{defPhi} and the continuous
            embedding $\left(\D,\|\cdot\|_{\D}\right)
            \hookrightarrow \left(L^{2^*},|\cdot|_{2^*}\right)$, for a suitable
            $k>0$ we have
                \begin{eqnarray*}
                    \widehat J(u)  &  =  &\frac{1}{2}\intrn|\nabla
                    u|^2\,dx-\frac{1}{2}\intrn f\big(|u+\nabla\Phi(u)|^2\big)\,dx\\
                                   & \geq&\frac{1}{2}\intrn|\nabla
                    u|^2\,dx-\frac{1}{2}\intrn f(|u|^2)\,dx\\
                                   & \geq&\frac{1}{2}\intrn|\nabla
                    u|^2\,dx-\frac{c_2}{2}\int_{\Omega_u}|u|^p\,dx-\frac{c_2}{2}
                    \int_{\R^n-\Omega_u}|u|^q\,dx\\
                                   & \geq&\frac{1}{2}\intrn|\nabla
                    u|^2\,dx-\frac{c_2}{2}\int_{\Omega_u}|u|^{2^*}\,dx-\frac{c_2}{2}
                    \int_{\R^n-\Omega_u}|u|^{2^*}\,dx\\
                                   &  =  &\frac{1}{2}\|u\|^2_{\D}-
                    \frac{c_2}{2}|u|^{2^*}_{2^*}\geq
                    \frac{1}{2}\|u\|_{\D}^2-k\|u\|_{\D}^{2^*}.
                \end{eqnarray*}
            Then $\widehat J(u)>C$ for $u\in S_\rho$ with $\rho$
            small enough.
        \end{proof}
        Now, before we prove that also the second geometrical
        assumption of the symmetrical mountain pass
        theorem is satisfied, we need a preliminary result.\\
        For every $\gamma>1$ and $u\in\V$ set
            \begin{equation}\label{deftilF}
                \widetilde F_u:w\in\W\mapsto\|u+\nabla
                w\|^\gamma_{L^p+L^q}.
            \end{equation}
        We have the following
        \begin{lemma}\label{deftilPhi}
            For every $u\in\V$ there exists a unique $\Phi_\gamma(u)\in\W$ s.t.
                \begin{equation}
                    \widetilde
                    F_u\big(\Phi_\gamma(u)\big)=\min_{w\in\W}\widetilde
                    F_u(w).
                \end{equation}
            Moreover, for every $V\subset\V$ s.t. $\dim V<+\infty$ we
            have
                \begin{equation}\label{propertyonV}
                    \exists\,\widetilde C_\gamma(V)>0\hbox{ s.t.
                    }\|u+\nabla\Phi_\gamma(u)\|^\gamma_{L^p+L^q}\geq
                    \widetilde C_\gamma\|u\|^\gamma_{\D}
                \end{equation}
            uniformly for $u\in V.$
        \end{lemma}
        \begin{proof}
            Since $\widetilde F_u$ is strictly convex, continuous
            and coercive on $\W$, by Weierstrass theorem there exists a
            unique minimizer $\Phi_\gamma(u).$

            Actually the minimizing map
                \begin{equation*}
                    \Phi_\gamma:u\rightarrow\Phi_\gamma(u)
                \end{equation*}
            is compact from $\V$ into $\W$.\\
            In fact, consider
                \begin{equation}\label{wconvinV}
                    u_n\rightharpoonup u\hbox{ in }\V.
                \end{equation}
            Since $\V\hookrightarrow L^p+L^q$ compactly, certainly
                \begin{equation}\label{convinLp+Lq}
                    u_n\rightarrow u\hbox{ in }L^p+L^q.
                \end{equation}
            Moreover, by the definition of $\Phi_\gamma,$
                \begin{equation*}
                    0\leq\|u_n+\nabla\Phi_\gamma(u_n)\|^\gamma_{L^p+L^q}\leq
                    \|u_n\|^\gamma_{L^p+L^q}
                \end{equation*}
            so,
                \begin{equation}\label{boundedinLp+Lq}
                    u_n+\nabla\Phi_\gamma(u_n)\hbox{ is bounded in
                    }L^p+L^q.
                \end{equation}
            From \eqref{convinLp+Lq} and \eqref{boundedinLp+Lq} we
            deduce that $\{\Phi_\gamma(u_n)\}$ is bounded in $\W$,
            so there exists $\overline w\in\W$ s.t. (up to a
            subsequence)
                \begin{equation}\label{wconvinLp+Lq}
                    \nabla\Phi_\gamma(u_n)\rightharpoonup\nabla\overline
                    w\hbox{ in }L^p+L^q.
                \end{equation}
            Now we prove that
                \begin{enumerate}
                    \item
                    $\displaystyle\lim_n\|u_n+\nabla\Phi_\gamma(u_n)\|_{L^p+L^q}=
                    \|u+\nabla\Phi_\gamma(u)\|_{L^p+L^q}$;
                    \item $\nabla\Phi_\gamma(u_n)\rightharpoonup
                    \nabla\Phi_\gamma(u)$ in $L^p+L^q.$
                \end{enumerate}
            Observe that, by the definition of $\Phi_\gamma$ and
            the triangular inequality,
                \begin{eqnarray*}
                    \|u+\nabla\Phi_\gamma(u)\|^\gamma_{L^p+L^q}
                    &  \leq  &
                    \|u+\nabla\Phi_\gamma(u_n)\|^\gamma_{L^p+L^q}\\
                    &  \leq  &
                    \big(\|u-u_n\|_{L^p+L^q}+\|u_n+
                    \nabla\Phi_\gamma(u_n)\|_{L^p+L^q}\big)^\gamma
                \end{eqnarray*}
            and then, by \eqref{convinLp+Lq}
                \begin{equation}\label{contofnorm}
                    \|u+\nabla\Phi_\gamma(u)\|^\gamma_{L^p+L^q}\leq\liminf_n
                    \|u_n+\nabla\Phi_\gamma(u_n)\|^\gamma_{L^p+L^q}.
                \end{equation}
            On the other hand, by definition of $\Phi_\gamma$
                \begin{equation*}
                    \|u_n+\nabla\Phi_\gamma(u_n)\|^\gamma_{L^p+L^q}\leq
                    \|u_n+\nabla\Phi_\gamma(u)\|^\gamma_{L^p+L^q}
                \end{equation*}
            and then, by \eqref{convinLp+Lq}
                \begin{equation}\label{contofnorm2}
                    \limsup_n\|u_n+\nabla\Phi_\gamma(u_n)\|^\gamma_{L^p+L^q}\leq
                    \|u+\nabla\Phi_\gamma(u)\|^\gamma_{L^p+L^q}.
                \end{equation}
            The claim $1$ follows from \eqref{contofnorm} and
            \eqref{contofnorm2}.\\
            Since $\|\cdot\|^\gamma_{L^p+L^q}$ is weakly
            lower semicontinuous, from \eqref{convinLp+Lq}, \eqref{wconvinLp+Lq}
            and the claim $1$ we deduce
                \begin{equation}\label{claim2}
                    \|u+\nabla\overline w\|^\gamma_{L^p+L^q}\leq
                    \liminf_n\|u_n+\nabla\Phi_\gamma(u_n)\|^\gamma_{L^p+L^q}
                    =\|u+\nabla\Phi_\gamma(u)\|^\gamma_{L^p+L^q}.
                \end{equation}
            By the uniqueness of the minimizer of $\widetilde
            F_u$, the inequality \eqref{claim2}
            implies that $\overline w=\Phi_\gamma(u)$ and then the claim $2$ is a
            consequence of \eqref{wconvinLp+Lq}.

            By a well known theorem, the claims $1$ and  $2$ and \eqref{convinLp+Lq}
            imply that
                \begin{equation*}
                    \nabla\Phi_\gamma(u_n)\rightarrow\nabla\Phi_\gamma(u)
                    \hbox{ in }L^p+L^q
                \end{equation*}
            and then $\Phi_\gamma$ is compact.

            Now, let $V\subset\V$ s.t. $\dim V<+\infty.$ By
            Weierstrass theorem
                \begin{equation}\label{deftilC}
                    \exists \widetilde C_\gamma:=\min_{\|u\|_{\D}=1 \atop
                    u\in V}\|u+\nabla\Phi_\gamma(u)\|^\gamma_{L^p+L^q}\geq
                    0.
                \end{equation}
            Actually, $\widetilde C_\gamma>0.$ In fact, if $\widetilde
            C_\gamma=0$, then there should exist $\overline u\in V$ s.t.
            $\|\overline u\|_{\D}=1$ and $\overline
            u+\nabla\Phi_\gamma(\overline u)=0,$ but it is not
            possible as we have already seen in the proof of
            Theorem \ref{critJ}. Now, if we consider $u\in
            V-\{0\}$ and set $\tilde u=u/\|u\|_{\D}$, since
            $\|\tilde u\|_{\D}=1$, we have that
                \begin{equation}\label{secondpartoftheo}
                    \frac{\|u+\nabla\Phi_\gamma(u)\|^\gamma_{L^p+L^q}}
                    {\|u\|^\gamma_{\D}}=\Big\|\tilde u+\nabla\Big(
                    \textstyle\frac{\Phi_\gamma(u)}{\|u\|_{\D}}\Big)\Big\|^\gamma_{L^p+L^q}
                    \geq\|\tilde u+\nabla\Phi_\gamma(\tilde
                    u)\|^\gamma_{L^p+L^q}\geq\widetilde C_\gamma.
                \end{equation}
            So \eqref{propertyonV} follows from
            \eqref{secondpartoftheo}.
        \end{proof}
        \begin{teo}\label{secondgeomount}
            For all $V\subset \V$ s.t. $\dim V<+\infty:$
            $\sup_{u\in V}\widehat J(u)<+\infty.$
        \end{teo}
        \begin{proof}
            Let $V\subset \V$ s.t. $\dim V<+\infty.$ Consider
            $u\in V$ and set
                $$\Omega:=\{x\in\R^n|\,|(u+\nabla\Phi(u))(x)|>1\}.$$
            Since inequality \eqref{controlnorm} implies that
                \begin{eqnarray*}
                    \|u+\nabla\Phi(u)\|_{L^p+L^q}^p  &  \leq  &
                    |u+\nabla\Phi(u)|_{L^p(\Omega)}^p\\
                    \noalign{\hbox{or}}
                    \|u+\nabla\Phi(u)\|_{L^p+L^q}^q  &  \leq  &
                    |u+\nabla\Phi(u)|_{L^q(\R^n-\Omega)}^q\;,
                \end{eqnarray*}
            certainly
                \begin{eqnarray}\label{maxleqmin}
                    \min\Big(\|u \!\!\! &  +  &\!\!\!  \nabla\Phi(u)\|_{L^p+L^q}^p    ,
                    \|u+\nabla\Phi(u)\|_{L^p+L^q}^q\Big)\nonumber\\
                    \nonumber\\
                    &  \leq  &
                    \max \Big(|u+\nabla\Phi(u)|_{L^p(\Omega)}^p,
                    |u+\nabla\Phi(u)|_{L^q(\R^n-\Omega)}^q\Big).
                \end{eqnarray}
            By \eqref{maxleqmin} and Lemma \ref{deftilPhi}
                \begin{align*}
                    \intrn f(|u+\nabla\Phi(u)|^2)\,dx &  \geq
                    c_1\int_\Omega
                    |u+\nabla\Phi(u)|^p\,dx+c_1\int_{\R^n-\Omega}
                    |u+\nabla\Phi(u)|^q\,dx\\
                    &  =
                    c_1|u+\nabla\Phi(u)|^p_{L^p(\Omega)}+c_1
                    |u+\nabla\Phi(u)|^q_{L^q(\R^n-\Omega)}\\
                    & \geq c_1\max\Big(|u+\nabla\Phi(u)|_{L^p(\Omega)}^p,
                    |u+\nabla\Phi(u)|_{L^q(\R^n-\Omega)}^q\Big)\\
                    & \geq c_1\min\Big(\|u +\nabla\Phi(u)\|_{L^p+L^q}^p,
                    \|u+\nabla\Phi(u)\|_{L^p+L^q}^q\Big)\\
                    & \geq c_1\min\Big(\|u +\nabla\Phi_p(u)\|_{L^p+L^q}^p,
                    \|u+\nabla\Phi_q(u)\|_{L^p+L^q}^q\Big)\\
                    & \geq c_1\min\big(\widetilde
                    C_p\|u\|_{\D}^p,\widetilde
                    C_q\|u\|_{\D}^q\big)\\
                    & \geq c_1\min(\widetilde
                    C_p\,,\widetilde C_q)\min\big(\|u\|_{\D}^p,\|u\|_{\D}^q\big),
                \end{align*}
            and then
                \begin{align}\label{Jbounded}
                    \widehat J(u) &  =  \frac{1}{2}\|u\|_{\D}^2-\frac{1}{2}\intrn
                    f(|u+\nabla\Phi(u)|^2)\,dx\nonumber\\
                                  &  \leq\frac{1}{2}\|u\|_{\D}^2-
                    c_1\min(\widetilde
                    C_p\,,\widetilde C_q)\min(\|u\|_{\D}^p,\|u\|_{\D}^q).
                \end{align}
            Since $2<p<q,$ we get our
            conclusion from \eqref{Jbounded}.
        \end{proof}

        \begin{proof}[Proof of Theorem \ref{th:maintheorem}]
            Since $\widehat J$ is $C^1$ and even, by
            Theorem \ref{JsatP.S:}, Theorem \ref{firstgeomount},
            Theorem \ref{secondgeomount} and the
            symmetrical version of the mountain pass theorem (see
            \cite{A.R.,B.B.F.}) certainly it possesses infinitely many critical
            points. Then the conclusion is a consequence of
            Theorem \ref{critJ}.
        \end{proof}

{\bf Acknowledgement}\\
The author is very grateful to Donato Fortunato and Alessio
Pomponio for their precious help and support.

\section*{Appendix}
As we have seen, in order to have infinitely many solutions for
the problem \eqref{eq:myprob} we need some assumptions on the
growth and on the convexity of the nonlinearity.  Here we want to
show an example of function satisfying those assumptions.

Consider the function $f:[0,+\infty[\rightarrow\R$ s.t.
    \begin{equation}\label{deff}
        f(x)=\left\{
            \begin{array}{lcl}
                ax^p+b & \hbox{if} & x>1\\
                cx^q      & \hbox{if} & x\leq1
            \end{array}
        \right.
    \end{equation}
where $2<p<2^*<q$ and the set of three numbers
$(a,b,c)\in\R^2\times]0,+\infty[$ is any solution of the system
    \begin{equation}\label{abcd}
        \left\{
            \begin{array}{lcl}
                a+b & = & c\\
                ap  & = & cq
            \end{array}
        \right..
    \end{equation}

    \begin{lemA}\label{lemmaA}
        There exist $\delta>0$ and $K_1>0$ s.t. $\forall
        (x,y)\in\;]1,1+\delta[\;\times\;]1-\delta,1[$
            \begin{equation}\label{fromlimit}
                f(x)-f(y)-\big(f'(y)|x-y\big)\geq K_1|x-y|^q.
            \end{equation}
    \end{lemA}
    \begin{proof}
        Consider the function $h:\;]1,+\infty[\times]0,1]$ s.t.
            \begin{equation}\label{definitionh}
                h(x,y)=\frac{f(x)-f(y)-\big(f'(y)|x-y\big)}{|x-y|^q}.
            \end{equation}
        that is
            \begin{equation*}
                h(x,y)=\frac{ax^p+b+(q-1)cy^q-qcxy^{q-1}}{|x-y|^q}.
            \end{equation*}
        Dividing numerator and denominator by $y^q$ and setting $z=x/y$, we
        get the new function
            \begin{equation}
                \tilde h(z,y)=\frac{az^py^{p-q}+by^{-q}+(q-1)c-qcz}{|z-1|^q}
            \end{equation}
        defined in the domain
        $\big\{(z,y)\in\;]1,+\infty[\times]0,1]\;\mid y>1/z\big\}.$\\
        We claim that
            \begin{equation*}
                \forall z>1:\quad \tilde h(z,1)=\min_{y> 1/z} \tilde
                h(z,\cdot).
            \end{equation*}
        We compute
            \begin{eqnarray}\label{parth}
                \frac{\partial \tilde h}{\partial y}(z,y)
                & = &
                \frac{a(p-q)z^py^{p-q-1}-bqy^{-q-1}}{|z-1|^q}\nonumber\\
                & = &
                \frac{a(p-q)z^py^{p}-bq}{y^{q+1}|z-1|^q}
                =\frac{g(zy)}{y^{q+1}|z-1|^q}
            \end{eqnarray}
        where $g(t) = a(p-q)t^p-bq.$\\
        By \eqref{abcd} we deduce that
            \begin{eqnarray*}
                g(1)   &  =  &  0\\
                g'(t)  &  <  &  0\quad\hbox{if }t>1
            \end{eqnarray*}
        so $g(zy)< 0$ because $zy>1.$ By \eqref{parth} we
        can conclude that the function $\tilde h (z,\cdot)$ is
        decreasing in $\left]\frac{1}{z},1\right]$ and then
            \begin{equation}\label{minh}
                \forall z>1:\quad\tilde h(z,y)\geq\tilde h(z,1).
            \end{equation}
        Now, by \eqref{minh} and using twice De l'H$\hat{\hbox{o}}$pital's rule,
        we compute
            \begin{eqnarray*}
                \lim_{(x,y)\rightarrow(1^+,1^-)}h(x,y)
                &  =  &
                \lim_{(z,y)\rightarrow(1^+,1^-)}\tilde h(z,y)\\
                &\geq &
                \lim_{(z,y)\rightarrow(1^+,1^-)}\tilde h(z,1)\\
                &  =  &
                \lim_{z\rightarrow1^+}\frac{az^p+b+(q-1)c-qcz}{(z-1)^q}\\
                &  =  &
                \lim_{z\rightarrow1^+}\frac{ap(p-1)z^{p-2}}{q(q-1)(z-1)^{q-2}}=+\infty.
            \end{eqnarray*}
    The inequality \eqref{fromlimit} is a consequence of the
    previous limit.
    \end{proof}

    \begin{teoA}
        There exists $K_2>0$ s.t. for every nonnegative numbers
        $x,y$
            \begin{equation}\label{ineqforf}
                f(x)-f(y)-f'(y)(x-y)\geq K_2 \min
                (|x-y|^p,|x-y|^q).
            \end{equation}
    \end{teoA}
    \begin{proof}
        We distinguish the following three cases
            \begin{enumerate}
                \item $0\leq y\leq 1<x$ or $0\leq x\leq 1<y$;
                \item $1<x,y$;
                \item $0\leq x,y\leq 1$.
            \end{enumerate}
            \begin{enumerate}
                \item If $0\leq y\leq 1<x$,
                then we consider these three possibilities
                    \begin{itemize}
                        \item $(x,y)    \in  \;  ]1,1+\delta[\;\times\;]1-\delta,1[$;
                        \item $(x,y)    \in  \;
                        ]1,1+\delta]\;\times\;[0,1-\delta]$
                        \item $(x,y)    \in  \;
                        [1+\delta,+\infty[\;\times\;[0,1],$
                    \end{itemize}
                where $\delta$ is the same as in Lemma A.\ref{lemmaA}.

                By Lemma A.\ref{lemmaA}, certainly \eqref{ineqforf} holds
                in $]1,1+\delta[\;\times\;]1-\delta,1[.$\\
                Since
                the function
                $h$ defined in \eqref{definitionh} is continuous in
                $[1,1+\delta]\;\times\;[0,1-\delta]$, by
                Weierstrass' theorem
                    \begin{equation*}
                        \exists
                        \min\{h(x,y)\mid(x,y)\in\;[1,1+\delta]\;\times\;[0,1-\delta]\}
                    \end{equation*}
                and then the inequality \eqref{ineqforf} holds
                also in $]1,1+\delta]\;\times\;[0,1-\delta].$\\
                Finally, suppose
                $(x,y)\in\;[1+\delta,+\infty[\;\times\;[0,1].$ Since for
                every $x\in[1+\delta,+\infty[$
                    \begin{equation*}
                        \min_{y\in[0,1]}y^{q-1}(c(q-1)y-cqx)=c(q-1)-cqx,
                    \end{equation*}
                then, by \eqref{abcd},
                    \begin{equation}\label{extimate}
                        \begin{array}{lcl}
                            f(x)-f(y)-f'(y)(x-y)  \!\!\!&  =  &
                            \!\!\!ax^p+b\\
                            &&\phantom{ax^p}
                            \!\!\!+y^{q-1}(c(q-1)y-cqx)\\
                            \!\!\!& \geq &\!\!\! ax^p+b+c(q-1)-cqx\\
                            \!\!\!& = &\!\!\! ax^p-ap(x-1)-a.
                        \end{array}
                    \end{equation}
                But
                    \begin{equation}\label{inf>0}
                        C_1:=\inf_{x\geq1+\delta}\frac{ax^p-ap(x-1)-a}{x^p}>0
                    \end{equation}
                so, by \eqref{extimate} and \eqref{inf>0},
                    \begin{equation*}
                        f(x)-f(y)-f'(y)(x-y)\geq C_1x^p\geq
                        C_1(x-y)^p
                    \end{equation*}
                and then the inequality \eqref{ineqforf} holds
                also in $[1+\delta,+\infty[\;\times\;[0,1].$

                We can use similar arguments for the case $0\leq x\leq
                1<y$.
                \item Suppose $1<x,y$. We have that
                    \begin{equation}\label{secondcase}
                        f(x)-f(y)-f'(y)(x-y)=a(x^p-y^p-py^{p-1}(x-y)).
                    \end{equation}
                In \cite{Kan} (see the proof of theorem 4, Chapter VIII) the following
                inequality has been proved: for all $r>2$ there exists a
                positive constant $C_2(r)$ s.t. for any $u\in\R$
                    \begin{equation}\label{Kantineq}
                        |u+1|^r\geq 1+ru+C_2(r)|u|^r.
                    \end{equation}
                If we set $r=p$ and replace $u$
                 by $\frac{x-y}{y}$, then by some calculus we get
                    \begin{equation}\label{Kantineq1}
                        x^p\geq
                        y^p+py^{p-1}(x-y)+C_2(p)|x-y|^p.
                    \end{equation}
                Inequality \eqref{ineqforf} follows from \eqref{secondcase}
                and \eqref{Kantineq1}.
                \item Suppose $0\leq x,y\leq 1$. Then
                    \begin{equation}
                        f(x)-f(y)-f'(y)(x-y)=c(x^q-y^q-qy^{q-1}(x-y))
                    \end{equation}
                so we get again \eqref{ineqforf} using
                \eqref{Kantineq} as before.
            \end{enumerate}
    \end{proof}
    \begin{teoA}
        Let $\hat f$ be the even extension of $f$, i.e.
            \begin{equation}
                \hat f(x)=\left\{
                    \begin{array}{ccl}
                        f(x) & \hbox{if}  & x\geq0\\
                        f(-x)& \hbox{if}  & x<0
                    \end{array}
                .\right.
            \end{equation}
        Then there exists $K_3>0$ s.t. for all $(x,y)\in\R^2$
            \begin{equation}\label{ineqforhatf}
                \hat f(x)-\hat f(y)-\hat f'(y)(x-y)\geq K_3 \min
                (|x-y|^p,|x-y|^q).
            \end{equation}
    \end{teoA}
    \begin{proof}
        We distinguish some cases
            \begin{itemize}
                \item $x,y\leq0$.

                Since $\hat f$ is even, certainly $\hat f'$ is odd and then,
                by \eqref{ineqforf},
                    \begin{eqnarray*}
                        \hat f(x)-\hat f(y) -\hat f'(y)(x-y)
                        &  =  &
                        \hat f(-x)-\hat f(-y)-\hat
                        f'(-y)(-x-(-y))\\
                        &  =  &
                        f(-x)-f(-y)-f'(-y)(-x-(-y))\\
                        &\geq &
                        K_2\min (|-x-(-y)|^p,|-x-(-y)|^q)\\
                        &  =  &
                        K_2\min (|x-y|^p,|x-y|^q).
                    \end{eqnarray*}
                \item $x\leq 0$ and $y\geq0.$

                We have that
                    \begin{equation}\label{explicit}
                        \hat f(x)-\hat f(y)-\hat f'(y)(x-y)=\hat
                        f(x)-\hat f'(y)x-\hat f(y)+\hat f'(y)y.
                    \end{equation}
                Since $\hat f'(y)\geq0$, the property $f_3$ (that
                can be easily proved)
                implies that
                    \begin{equation}\label{firstinequality}
                        \hat f(x)-\hat f'(y)x\geq \hat f(x)\geq
                        c_1\min (|x|^p,|x|^q),
                    \end{equation}
                and, on the other hand, by \eqref{ineqforf}
                    \begin{equation}\label{secondinequality}
                        \begin{array}{lcl}
                            -\hat f(y)+\hat
                            f'(y)y
                            &  =  &
                            f(0)-f(y)-f'(y)(0-y)\\
                            \\
                            &\geq &
                            K_2\min(|y|^p,|y|^q).
                        \end{array}
                    \end{equation}
                Comparing \eqref{explicit},
                \eqref{firstinequality} and
                \eqref{secondinequality} we get
                    \begin{eqnarray*}
                        \hat f(x)-\hat f(y) -\hat f'(y)(x-y)
                        & \geq  &
                        C_3\big(\min(|x|^p,|x|^q)+\min(|y|^p,|y|^q)\big)\\
                        & \geq  &
                        C_4\min(|x|^p+|y|^p,|x|^q+|y|^q)\\
                        & \geq  &
                        C_5\min\big((|x|+|y|)^p,(|x|+|y|)^q\big)\\
                        &   =   &
                        C_5\min(|x-y|^p,|x-y|^q).
                    \end{eqnarray*}
                where $C_3,$ $C_4$ and $C_5$ are positive
                constants.
                \item $x\geq0$ and $y\leq0$.

                The inequality \eqref{ineqforhatf} can be proved
                by similar arguments as before.
                \item $x,y\geq0$.

                The inequality \eqref{ineqforhatf} follows
                directly from \eqref{ineqforf}.
            \end{itemize}
    \end{proof}
Finally, define $\overline f:\R^n\rightarrow\R$ as the radial
extension of $f$, namely
    \begin{equation}\label{defoverf}
        \overline f(x)=f(|x|),\quad\forall x\in\R^n.
    \end{equation}

    \begin{teoA}\label{teoA1}
        The function $\overline f$ defined by \eqref{defoverf}
        and \eqref{deff} satisfies the inequality
            \begin{equation}\label{ap:f_2}
                \overline f(x)-\overline f(y)-\big( \overline
                f\,'(y)\big|
                x-y\big)
                  =
                c_1\min\big(|x-y|^p,|x-y|^q\big)
            \end{equation}
        for some positive constant $c_1$ which doesn't depend on
        $x,y\in\R^n.$
    \end{teoA}

    \begin{proof}
        It is very easy to verify
        that $\overline f$ satisfies the inequality
            \begin{equation}\label{ap:ineq}
                 f(x)\geq c_1
                        \min(|x|^p,|x|^q)
            \end{equation}
        for all $x\in\R^n.$

        If $y=0$, then
        \eqref{ap:f_2} follows trivially from \eqref{ap:ineq}.

        If $y\neq 0$, observe that for all
        $x\in\R^n$
            \begin{equation}\label{explicit2}
                \begin{array}{lcl}
                    \overline f(x)-\overline f(y)-\big( \overline
                    f\,'(y)\!\!\!\!\!\!&\big|&\!\!\!\!\!\!
                    x-y\big)\\
                    \!\!\!&=&\!\!\!
                    f(|x|)-f(|y|)
                    \displaystyle-\frac{f'(|y|)}{|y|}(x|y)+f'(|y|)|y|.
                \end{array}
            \end{equation}
        Now consider the following three cases
            \begin{itemize}
                \item $x=ty,$ $t\geq0$;
                \item $x=ty,$ $t<0$;
                \item $x\neq ty,$ $t\in\R$.
            \end{itemize}
        If $x=ty$ for $t\geq0$, then $(x|y)=|x||y|$ and
        $|x-y|=\big||x|-|y|\big|$ so, by \eqref{explicit2} and \eqref{ineqforf},
            \begin{eqnarray*}
                \overline f(x)-\overline f(y)-\big( \overline
                f\,'(y)\big|
                x-y\big)
                &  =  &
                f(|x|)-f(|y|)-f'(|y|)(|x|-|y|)\\
                &\geq &
                K_2
                \min\big(\big||x|-|y|\big|^p,\big||x|-|y|\big|^q\big)\\
                &  =  &
                K_2
                \min(|x-y|^p,|x-y|^q).
            \end{eqnarray*}
        If $x=ty$ for $t<0$, then $(x|y)=-|x||y|$ and
        $|x-y|=|x|+|y|$ so, by \eqref{explicit2} and
        \eqref{ineqforhatf},
            \begin{eqnarray*}
                \overline f(x)-\overline f(y)-\big( \overline
                f\,'(y)\big|
                x-y\big)
                &  =  &
                \hat f(|x|)-\hat f(-|y|)-\hat
                f\,'(-|y|)(|x|-(-|y|))\\
                &\geq &
                K_3\min\big(\big||x|+|y|\big|^p,\big||x|+|y|\big|^q\big)\\
                &  =  &
                K_3\min\big(|x-y|^p,|x-y|^q\big).
            \end{eqnarray*}
        Finally, if $x\not\in\{ty|t\in\R\}$, then $x=x_1+x_2$
        where $x_1\in\{ty|t\in\R\}$ and $(x_2|y)=0$. Since $x_1\parallel y$,
        from the previous cases we have
            \begin{equation}\label{tilf2forx1}
                \overline f(x_1)-\overline f(y)-
                \big(\overline f\,'(y)\big|x_1-y\big)\geq C_6
                \min(|x_1-y|^p,|x_1-y|^q)
            \end{equation}
        where $C_6=\min(K_2,K_3).$
        Moreover, observe that
        for all $a,b\geq0$ the following inequality holds
            \begin{equation}\label{littleprop}
                f(\sqrt{a+b})\geq f(\sqrt a)+ f(\sqrt b),
            \end{equation}
        so, by \eqref{littleprop}, \eqref{tilf2forx1} and property
        $f_3$, we have
            \begin{eqnarray*}
                \overline f(x)  -  \overline f(y)&-&
                \big(\overline f\,'(y)\big|x-y\big)\\
                &  =  &
                f(\sqrt{|x_1|^2+|x_2|^2})-\overline f(y)-
                \big(\overline f\,'(y)\big|x_1-y\big)\\
                &\geq &
                f(|x_1|)+f(|x_2|)-\overline f(y)-
                \big(\overline f\,'(y)\big|x_1-y\big)\\
                &  =  &
                \overline f(x_1)-\overline f(y)-
                \big(\overline f\,'(y)\big|x_1-y\big)+\overline
                f(x_2)\\
                &\geq &
                C_6
                \min(|x_1-y|^p,|x_1-y|^q)+c_1\min(|x_2|^p,|x_2|^q)\\
                &\geq &
                C_7\min\big((|x_1-y|^2)^\frac{p}{2}+(|x_2|^2)^\frac{p}{2},
                (|x_1-y|^2)^\frac{q}{2}+(|x_2|^2)^\frac{q}{2}\big)\\
                &\geq &
                C_8\min\big((|x_1-y|^2+|x_2|^2)^\frac{p}{2},
                (|x_1-y|^2+|x_2|^2)^\frac{q}{2}\big)\\
                &  =  &
                C_8\min(|x-y|^p,|x-y|^q)
            \end{eqnarray*}
        where $C_7$ and $C_8$ are suitable positive constants.

    \end{proof}

    Now, consider the function $f:\R\rightarrow\R$ s.t.
    \begin{equation*}
        f(x)=\left\{
            \begin{array}{lcl}
                a|x|^\frac{p}{2}+b & \hbox{if} & |x|>1\\
                c|x|^\frac{q}{2}   & \hbox{if} & |x|\leq1
            \end{array}
        \right.
    \end{equation*}
where $2<p<2^*<q$ and $(a,b,c)\in\R^2\times]0,+\infty[$ is any
solution of the system
    \begin{equation*}
        \left\{
            \begin{array}{lcl}
                a+b & = & c\\
                ap  & = & cq
            \end{array}
        \right..
    \end{equation*}
It is easy to verify that $f$ satisfies $f_1,$ $f_3$ and $f_4$.\\
Moreover, applying Theorem A.\ref{teoA1} to the function
    $$g:\xi\in\R^n\mapsto f((\xi|\xi))\in\R,$$
we verify that $f$ satisfies also $f_2.$

\end{document}